\documentclass[preprint,11pt]{elsarticle}
\makeatletter
\def\ps@pprintTitle{%
 \let\@oddhead\@empty
 \let\@evenhead\@empty
 \def\@oddfoot{\centerline{\thepage}}%
 \let\@evenfoot\@oddfoot}
\makeatother
\usepackage{amssymb}
\usepackage{amsmath}
\usepackage{amsthm}
\usepackage{epsfig}
\usepackage{color}
\usepackage{makeidx}
\usepackage{bbm}
\usepackage[top=1in, bottom=1in, inner=1in, outer=1in]{geometry}
\usepackage{setspace}

\newtheorem{theorem}{Theorem}

\newcommand*\xbar[1]{%
  \hbox{%
    \vbox{%
      \hrule height 0.5pt 
      \kern0.4ex
      \hbox{%
        \kern-0.15em
        \ensuremath{#1}%
        \kern-0.15em
      }%
    }%
  }%
}

\begin{document}

\begin{frontmatter}

\title{Subdiffusivity of Brownian Motion among a Poissonian Field of Moving Traps}

\author{Mehmet \"{O}z}
\ead{mehmet.oz@ozyegin.edu.tr}

\address{Department of Natural and Mathematical Sciences, Faculty of Engineering, \"{O}zye\u{g}in University, Istanbul, Turkey}

\begin{abstract}
Our model consists of a Brownian particle $X$ moving in $\mathbb{R}$, where a Poissonian field of moving traps is present. Each trap is a ball with constant radius, centered at a trap point, and each trap point moves under a Brownian motion independently of others and of the motion of $X$. Here, we investigate the `speed' of $X$ on the time interval $[0,t]$ and on `microscopic' time scales given that $X$ avoids the trap field up to time $t$. Firstly, following the earlier work of Athreya et al.\ [{\em Math.\ Phys.\ Anal.\ Geom.}\ \textbf{20}:1 (2017)], we obtain bounds on the maximal displacement of $X$ from the origin. Our upper bound is an improvement of the corresponding bound therein. Then, we prove a result showing how the speed on microscopic time scales affect the overall macroscopic subdiffusivity on $[0,t]$. Finally, we show that $X$ moves subdiffusively even on certain microscopic time scales, in the bulk of $[0,t]$. The results are stated so that each gives an `optimal survival strategy' for the system. We conclude by giving several related open problems.
\end{abstract}

\vspace{3mm}

\begin{keyword}
Random motion in random environment \sep Poissonian traps \sep moving trap field  \sep subdiffusive \sep optimal survival strategy. 
\vspace{3mm}
\MSC[2010] 60K37 \sep 60D05 \sep 60K35  \sep 82C22 
\end{keyword}

\end{frontmatter}

\pagestyle{myheadings}
\markright{SUBDIFFUSIVITY OF BM IN MOVING TRAP FIELD\hfill}

\section{Introduction}\label{intro}

Brownian motion among a moving Poissonian trap field attached to $\mathbb{R}^d$ has been studied recently in \cite{P2012}. The discrete analogue of this model, that is, random walk among a Poisson system of moving traps attached to the discrete lattice $\mathbb{Z}^d$, has been studied concurrently in \cite{R2012}. Both works gave results on the large time asymptotics of the survival probability of the randomly moving particle, where one defines survival up to time $t$ to be the event that the particle has not hit the traps until that time. Another important problem on this model is that of optimal survival: How must have the system, which is composed of the randomly moving particle and the field of traps, behaved (what strategy must it have followed) given that the particle has avoided the traps up to time $t$? In the discrete setting, Athreya et al. addressed this issue in \cite{ADS2016}, where they considered the one-dimensional model on $\mathbb{Z}$ and focused on the maximal displacement of the random walk from origin. Conditioned on survival, in order to avoid the traps, it is natural to expect the particle to go not as far from origin as it otherwise would in the absence of traps. In \cite{ADS2016}, Athreya et al. indeed showed that for large $t$, with overwhelming probability, the random walk behaves subdiffusively in order to avoid traps. In the present work, we extend their result to the continuous case of Brownian motion among a Poissonian trap field on $\mathbb{R}$, with an improvement on the upper bound of maximal displacement from origin. Moreover, we show that if the particle is `too slow' or `too fast' when traversing distances on the order of its maximal displacement, its overall subdiffusivity is strengthened; therefore, in particular the way in which $X$ reaches the point of its maximal displacement also matters. Finally, we consider microscopic time scales and show that the Brownian particle must behave subdiffusively even on certain time scales of order $o(t)$.
    
\subsection{Formulation of the problem}		
		
The setting of a Brownian particle among a moving random trap field is formed as follows. Let $X$ be a Brownian particle and let $X=(X_s)_{s\geq 0}$ represent its path on $\mathbb{R}$, where we take $X_0=0$. Note that we use $X$ both as the name of the particle and as the random variable representing its sample path. Let $\mathbb{P}^X$ and $\mathbb{E}^X$ denote respectively the probability law and corresponding expectation for $X$. Create a random environment on $\mathbb{R}$ via a `dynamic' Poisson point process (PPP) $\Pi=(\Pi_s)_{s\geq 0}$ as follows. Let $\Pi_0=\left\{x_i\right\}_i$ be a PPP with constant intensity $\lambda>0$, placed on $\mathbb{R}$ at time $t=0$. Here, we refer to the points of the PPP as trap points. Now let each trap point $x_i$ at $t=0$ move under a Brownian motion $Y^i=(Y^i_s)_{s\geq 0}$ independently of all others and of $X$ so that $\Pi_s:=\left\{x_i+Y^i_s\right\}_i$ is the point process after each $x_i$ has moved for time $s$. Applying the mapping theorem for Poisson processes, it can be shown that for each $s>0$, $\Pi_s$ is also a PPP with the initial intensity $\lambda$ (see \cite{BMW1997} for details). By a `trap' associated to a trap point at $x\in\mathbb{R}$, we mean a closed ball of fixed radius $a>0$ centered at $x$ (in $d=1$, note that this is just a closed interval). Then, the moving (random) trap field $K=(K_s)_{s\geq 0}$ is given by 
\begin{equation} K_s:=\underset{x_i \in\: \text{supp}(\Pi_0)}{\bigcup} \bar{B}(x_i+Y^i_s,a), \nonumber
\end{equation}     
where $B(x,a)$ denotes the open ball centered at $x$ with radius $a$, and $\bar{A}$ denotes the closure of a set $A\subseteq\mathbb{R}$. Let $\mathbb{P}$ and $\mathbb{E}$ denote respectively the probability law and corresponding expectation for the moving trap field, that is, for $\Pi=(\Pi_s)_{s\geq 0}$.

In this work, the Brownian particle $X$ is assumed to live in $\mathbb{R}$ with the trap field $K$ attached to it. Let $R(X)=(R_t(X))_{t\geq 0}$ represent the range process of $X$, where $R_t(X):=\left\{X_s:0\leq s \leq t\right\}$ is the set of all points in $\mathbb{R}$ that $X$ visits up to time $t$. Define $T=\inf\left\{s\geq 0:R_s(X)\cap K_s\neq\emptyset\right\}$ to be the first time that $X$ hits a trap, and $\left\{T>t\right\}$ to be the event of survival up to time $t$.

The probability measure of interest is $(\mathbb{E}\times \mathbb{P}^X)(\:\cdot \mid T>t)$, the annealed probability conditioned on survival of $X$ up to time $t$. By an \textit{optimal survival strategy}, we mean a collection of events $\left\{A_t\right\}_{t>0}$ indexed by $t$ such that 
\begin{equation} \underset{t\rightarrow\infty}{\lim}(\mathbb{E}\times \mathbb{P}^X)\left(A_t \mid T>t\right)=1. \nonumber
\end{equation}
We look for optimal survival strategies concerning the `speed' of $X$, where by speed, we refer to the Lebesgue measure of the range per unit time over a given time interval. Our main result, Theorem~\ref{teo1}, is on the maximal displacement of $X$ from origin up to time $t$, similar to \cite[Thm.1.2]{ADS2016}. We emphasize that this is a result on a `macroscopic' time scale as it gives a strategy over the entire interval $[0,t]$. Theorem~\ref{teo3} studies the effect of a certain `microscopic' speed on the maximal displacement over $[0,t]$, which is a measure of the macroscopic speed. Lastly, in Theorem~\ref{teo2}, we show that trap-avoiding forces $X$ to move subdiffusively on certain microscopic time scales as well. 

\subsection{History}

Trapping problems in the context of a single randomly moving particle among a (Poissonian) random field of traps have a long history. In the continuous setting of a Brownian particle among frozen (static) Poissonian traps in $\mathbb{R}^d$, the large time survival asymptotics were studied in \cite{DV1975} and \cite{S1998}, and optimal survival strategies were studied in \cite{S1990}, \cite{S1991} and \cite{P1999} in dimensions $d=1$, $d=2$ and $d\geq 3$, respectively. Survival asymptotics involving a moving (dynamic) field of traps were studied in \cite{P2012}; however, as far as we know, there is no corresponding work on optimal survival strategies. In the discrete setting, where the continuum $\mathbb{R}^d$ is replaced by the integer lattice $\mathbb{Z}^d$, the survival asymptotics of a random walk among frozen Bernoulli traps were studied in \cite{DV1979} and \cite{A1995}. We note that in the discrete setting, as long as hard-killing rule is applied, where the system is killed instantly the first time it hits a trap, there is no difference between the cases of Bernoulli traps and Poissonian traps provided that the traps are frozen.

In the discrete setting, a dynamic version of the model studied in \cite{A1995} was introduced in \cite{R2012} as follows. At $t=0$, the random walker $X$ is placed at the origin, and the environment is composed of a trap field on $\mathbb{Z}$ with each integer site having a random number of trap points, which are i.i.d.\ with a Poisson distribution. The dynamics of each trap point and $X$ are governed by independent random walks, where trap points have a common jump rate, and $X$ has a different jump rate. In \cite{R2012}, both the annealed and quenched survival asymptotics were studied under a soft-killing rule, where the particle $X$ is killed at a rate proportional to the number of trap points present at the site visited and to an interaction parameter $\gamma\geq0$. (Note that by taking $\gamma=\infty$, one may switch to hard-killing rule.) In \cite{ADS2016}, Athreya et al. studied the optimal survival problem on the model introduced in \cite{R2012}, and found bounds on maximal displacement of $X$ from origin conditioned on survival up to time $t$. The current work originated from \cite{ADS2016} in search for an extension of the results therein to the continuous setting, and for a better upper bound on the maximal displacement of $X$.    
 
\vspace{0.3cm}

The organization of the paper is as follows. In Section~\ref{s2}, we present our results. In Section~\ref{s3}, we develop the preparation needed for the proofs of the results. Section~\ref{s4} is devoted to the proof of Theorem~\ref{teo1}, which is our main result. The proofs of Theorem~\ref{teo3} and Theorem~\ref{teo2} are given in Section~\ref{s5}. Finally, we briefly mention some open problems related to our model in Section~\ref{s6}.

\section{Results}\label{s2}

We introduce further notation in order to state our results. Let
\begin{equation} M_t(X):=\sup_{s\leq t} X_s, \quad m_t(X):=\inf_{s\leq t} X_s, \nonumber
\end{equation}
and define $\lVert X \rVert_t:=\max\left\{M_t(X),-m_t(X)\right\}$ to be the maximal displacement of $X$ from origin up to time $t$. We write $P:=\mathbb{E}\times\mathbb{P}^X$ as the annealed probability measure for ease of notation. 

We now present our main result, which identifies an optimal survival strategy for $X$ among the moving field of traps described in the introduction. The result below concerns the macroscopic behavior of $X$ conditioned on survival up to time $t$.

 \begin{theorem}\label{teo1} 
There exist constants $c_1>0$ and $c_2>0$ such that 
\begin{equation} \underset{t\rightarrow\infty}{\lim}P\left(\lVert X \rVert_t \in (c_1\,t^{1/3},c_2\,t^{5/11})\mid T>t\right)= 1. \label{eq01}
\end{equation} 
\end{theorem}
We emphasize that Theorem~\ref{teo1} is an extension of and improvement on the corresponding result in \cite{ADS2016}. It extends the corresponding result in \cite{ADS2016} to continuous setting, and it is an improvement in that the epsilon in the exponent of $t$ in the upper bound therein is lost, and that the exponent of $t$ in the upper bound is improved from $11/24$ to $5/11$. 

The next result addresses the following question: How does $X$ behave on its way to the point of maximal displacement from origin, given that it avoids the trap field up to time $t$? As the following theorem shows, if $X$ moves `too fast' or `too slow' when traversing distances on the order of its maximal displacement from origin, then its overall subdiffusivity on $[0,t]$ is strengthened.     

\begin{theorem}\label{teo3} 
Let $0<\varepsilon\leq 1$ and $0<\kappa<1$. Let $A_t$ be the event that $X$ traverses the spatial interval $(\kappa\lVert X \rVert_t,\lVert X \rVert_t)$ at least diffusively fast or stays inside this interval for a total time of length $\varepsilon t$. Then there exists a constant $c_3>0$ such that  
\begin{equation} \underset{t\rightarrow\infty}{\lim}P\left(A_t, \lVert X \rVert_t \geq c_3\,t^{4/9} \mid T>t \right)= 0 . \nonumber
\end{equation} 
\end{theorem}
We compare Theorem~\ref{teo3} to Theorem~\ref{teo1}, and see that the exponent in the upper bound for $\lVert X \rVert_t$ in Theorem~\ref{teo1} decreases from $5/11$ to $4/9$ in Theorem~\ref{teo3}, which means that traversing distances on the order of maximal displacement from origin diffusively (too fast in this case) or spending too much time far away from origin tends to confine $X$ to a smaller interval around origin, conditional on survival from traps. This could be heuristically explained as follows. When $X$ traverses large distances without being trapped, the trap points are swept out of the way (actually, the trap points move out of the way, away from $X$, at least as fast as $X$ moves) and piled up near the boundary of the range of $X$, but this sweeping away becomes probabilistically too costly if the distance traversed by $X$ diffusively is as large as $t^{4/9}$ up to a large enough constant. On the other hand, the trap points that are already piled up near the boundary of the range, while $X$ was on its way to maximal displacement, will catch up with $X$ if $X$ is too slow moving back towards origin. Therefore if $X$ spends too much time at distances on the order of its maximal displacement without being trapped, it can move at most $t^{4/9}$ away from origin up to a large enough constant. A larger displacement would mean piling up of too many trap points near the boundary of the range so that at least one catches up with $X$ with overwhelming probability.

The next result concerns the microscopic behavior of $X$ conditioned on survival up to time $t$. By `microscopic', we mean over time scales of order $o(t)$ as $t\rightarrow\infty$.  

\begin{theorem}\label{teo2} 
Let $\varepsilon>0$ and $f:\mathbb{R}_+\to\mathbb{R}_+$ be a function such that $f(t)\rightarrow\infty$ and $f(t)=o(t^{1/3})$ as $t\rightarrow\infty$. Let $B_t$ be the event that $X$ is at least diffusively fast on $\geq \varepsilon t^{1/3}/f(t)$ many pairwise disjoint intervals in $[0,t]$ of length $t^{2/3}f(t)$ each. Then,
\begin{equation} \underset{t\rightarrow\infty}{\lim}P\left(B_t\mid T>t\right)=0. \nonumber
\end{equation}  
\end{theorem}
Theorem~\ref{teo2} says that conditioned on survival up to time $t$, with overwhelming probability, $X$ is not diffusive on time scales of order higher than $t^{2/3}$ in the `bulk' of $[0,t]$ for large $t$. Hence, $X$ is subdiffusive not only on the macroscopic scale of $t$ but also on microscopic time scales as long as they are higher order than $t^{2/3}$.

\section{Preparations}\label{s3}

In this section, we aim at obtaining a suitable expression that will serve as an upper bound for 
\[P(X\in \cdot \mid T>t).  \]
The line of argument will be similar to the one in \cite{ADS2016}. For an upper bound, we write
\begin{equation} P(X\in \cdot \mid T>t)=\frac{P(T>t,\:X\in \cdot \:)}{P(T>t)}, \label{eq30}  
\end{equation}
rewrite the numerator by `integrating out' the Poisson field $\Pi$, and bound the denominator from below via a survival strategy that is not too costly.

Henceforth, $c$, $c_1$, $c_2$, etc. will denote generic constants, whose values may change from line to line. The notation $c(\kappa)$ will be used to mean that the constant $c$ depends on the parameter $\kappa$. Furthermore, we will use $B=(B_s)_{s\geq 0}$ to denote a generic standard Brownian motion, and $P_x$ to denote the law of $B$ started at position $x\in\mathbb{R}$. When random variables such as $R_t$, $M_t$, etc. are written without regard to a particular Brownian motion, they are to be understood as functions of $B$. We will use $\mathbbm{1}_E$ as the indicator function for an event $E$, and $|A|$ as the Lebesgue measure of a set $A\subset \mathbb{R}$.

The survival probability of Brownian motion among a Poissonian trap field is closely related to a particular functional of the Brownian path, namely the `Wiener sausage'. Let $Y=(Y_s)_{s\geq 0}$ be the path of a Brownian particle $Y$. Then the Wiener sausage associated to $Y$ up to time $t$ is defined as 
\[ W_0(t):=\underset{s\leq t}{\bigcup} B(Y_s,a). \]
If $f:[0,\infty)\to\mathbb{R}$ is a deterministic function, letting $f(s)=f_s$, the Wiener sausage associated to $Y$ \textit{with drift $f$} up to time $t$ can be defined as 
\[ W_f(t):=\underset{s\leq t}{\bigcup} B(Y_s+f_s,a). \]
Then, by a standard application of Fubini's theorem, one can integrate out the Poisson field and show as in \cite[Lemma 2.1]{P2012} that
\begin{equation} \mathbb{P}(T>t\mid X)=\exp\left(-\lambda\:\mathbb{E}^Y[\,|W_X(t)|\mid X]\right), \label{eq31}
\end{equation}
where $Y$ is independent of $X$, and $\mathbb{P}(\cdot \mid X)$ denotes the conditional probability given $X$. Thanks to (\ref{eq31}), instead of dealing with the entire trap field $\Pi$ and $X$, it is enough to deal with two independent Brownian motions $X$ and $Y$. Then, the numerator in (\ref{eq30}) can be written as 
\begin{align} P(T>t,\:X\in \cdot \:)=&\:(\mathbb{E}^X\times\mathbb{E})[\mathbbm{1}_{\left\{X\in\: \cdot \:\right\}} \mathbbm{1}_{\left\{T>t\right\}}] \nonumber \\
=&\: \mathbb{E}^X[\mathbbm{1}_{\left\{X\in\: \cdot \:\right\}} \mathbb{P}(T>t \mid X)] \nonumber \\
=&\: \mathbb{E}^X [\mathbbm{1}_{\left\{X\in\: \cdot \:\right\}} \exp\left(-\lambda \:\mathbb{E}^Y[\,|W_X(t)|\mid X] \right)]. \label{eq32}
\end{align}

Next, we obtain a lower bound for the denominator in (\ref{eq30}). Let $r=r(t)$ with $r(t)\rightarrow\infty$ as $t\rightarrow\infty$. One way for $X$ to survive is to be confined to the ball $B(0,r)$ while $B(0,r+a)$ stays free of trap points up to time $t$. Recall that $R_t(X)=\left\{X_s:0\leq s \leq t\right\}$.  By a standard result (see for example \cite{PS1978}) on the probability of confinement of a Brownian particle in a ball, for all $t>0$,
\begin{equation} P(R_t(X)\subseteq B(0,r))\geq c_d\  e^{-\frac{\rho_d}{r^2}t}, \label{eq33}
\end{equation}
where $\rho_d$ is the principal Dirichlet eigenvalue for the open unit ball in $\mathbb{R}^d$ and $c_d$ is a constant that depends on dimension. Since $d=1$ throughout this work, we suppress the dependence on dimension. Define $T_{B(0,r)}:=\inf\left\{s\geq 0:B(0,r)\cap K_s\neq\emptyset\right\}$ to be the first time that the trap field hits $B(0,r)$. By an extension of (\ref{eq31}), it is easy to see that (\cite[Remark 2.3]{P2012})
\begin{equation} \mathbb{P}(T_{B(0,r)}>t)=\exp\left(-\lambda\:\mathbb{E}^Y\left[\,\bigg|\underset{s\leq t}{\bigcup} B(Y_s,r+a)\bigg|\,\right]\right).  \label{eq34}
\end{equation}  
Now let $\varepsilon>0$ and write
\begin{equation} \underset{s\leq t}{\bigcup} B(Y_s,r+a)=r \underset{s\leq t}{\bigcup} B(Y_s/r,1+a/r) \label{eq35}
\end{equation}
so that for all large $t$, since $r(t)\rightarrow\infty$, we have 
\begin{equation} B(Y_s/r,1+a/r) \subset B(Y_s/r,1+\varepsilon) \label{eq36}
\end{equation}
for all $s$. Then, by Brownian scaling, it follows from (\ref{eq35}) and (\ref{eq36}) that for all large $t$,
\begin{align} \mathbb{E}^Y\left[\,\bigg|\underset{s\leq t}{\bigcup} B(Y_s,r+a)\bigg|\,\right]\leq &\:
r\, \mathbb{E}^Y\left[\,\bigg|\underset{u\leq t/r^2}{\bigcup} B(Y_u,1+\varepsilon)\bigg|\,\right]  \nonumber \\
=&\: r \left(\mathbb{E}^Y\left[\,|R_{t/r^2}(Y)|\,\right]+2(1+\varepsilon)\right)=\mathbb{E}^Y\left[\,|R_t(Y)|\,\right]+2r(1+\varepsilon), \label{eq37}
\end{align}
where the last equality follows from (1.4) in \cite{F1951}. Then, choosing $r(t)=t^{1/3}$ for optimality, it follows from (\ref{eq33}), (\ref{eq34}) and (\ref{eq37}) that there exists a constant $c(\lambda)$ such that for all large $t$,
\begin{equation} P(T>t)\geq e^{-c(\lambda) t^{1/3}} \exp\left\{-\lambda\:\mathbb{E}^Y\left[\,|R_t(Y)|\,\right]\right\}. \label{eq38}
\end{equation}
Finally, using (\ref{eq32}) and (\ref{eq38}), and noting that $|W_X(t)|=|R_t(Y+X)|+2a$, it follows from (\ref{eq30}) that
\begin{equation} P(X\in \cdot \mid T>t)\leq e^{c(\lambda) t^{1/3}}\mathbb{E}^X\left[\mathbbm{1}_{\left\{X\in\: \cdot\right\}}\exp\left\{-\lambda\:\mathbb{E}^Y\left[\,|R_t(Y+X)|-|R_t(Y)|\,\mid X\right] \right\}  \right] \label{eq1}
\end{equation}
for all large $t$. This will be the starting point in the proof of Theorem~\ref{teo1}, which is given in two parts in the next section.

\section{Proof of Theorem~\ref{teo1}}\label{s4}

\noindent \textbf{Proof of lower bound}

\noindent From \cite[Prop.4.4]{TV2006}, we have the following asymptotics for the range of $B$: for $a>0$,
\begin{equation} P_0(|R_t|<a)\sim 8\pi^2 \frac{t}{a^2} e^{-\frac{\pi^2}{2}\frac{t}{a^2}}, \nonumber
\end{equation}
where $\sim$ denotes asymptotic similarity as $t\rightarrow\infty$. Hence for any $c>0$, we have
\begin{equation} P_0(|R_t|<c t^{1/3})\sim \frac{8\pi^2}{c^2} t^{1/3} e^{-\frac{\pi^2}{2c^2}t^{1/3}}. \label{eq2}
\end{equation}
Now consider the second expectation on the right-hand side of (\ref{eq1}). Since both $X$ and $Y$ have continuous sample paths almost surely, with $\mathbb{P}^X$-probability $1$, 
\begin{align} \mathbb{E}^Y\left[\:|R_t(Y+X)|-|R_t(Y)|\:\right]&=\mathbb{E}^Y\left[\sup_{s\leq t} (Y_s+X_s)-\inf_{s\leq t} (Y_s+X_s)-\sup_{s\leq t} (Y_s)+\inf_{s\leq t} (Y_s)\right] \nonumber \\
&=\mathbb{E}^Y\left[\sup_{s\leq t} (Y_s+X_s)+\sup_{s\leq t} (-Y_s-X_s)-2\sup_{s\leq t} (Y_s)\right] \nonumber \\
&=\mathbb{E}^Y\left[\sup_{s\leq t} (Y_s+X_s)+\sup_{s\leq t} (Y_s-X_s)-2\sup_{s\leq t} (Y_s)\right] \geq 0 \label{eq3}
\end{align}
for all $t\geq 0$, where we have used the symmetry of $Y$ in passing to the second and third equalities, and concluded that the expression is non-negative since the first two suprema on the right-hand side could be opened at the argument maximum of $Y$ for a lower bound. 
Noting that $\mathbb{P}^X(\lVert X \rVert_t <ct^{1/3})\leq P_0(|R_t|< 2c t^{1/3})$, it follows from (\ref{eq1}), (\ref{eq2}) and (\ref{eq3}) that 
\begin{align} P(\lVert X \rVert_t <c t^{1/3} \mid T>t)&\leq e^{c(\lambda) t^{1/3}} \mathbb{P}^X(\lVert X \rVert_t <c t^{1/3}) \nonumber \\ 
&\leq e^{c(\lambda) t^{1/3}} e^{-\frac{\pi^2}{8 c^2}t^{1/3}(1+o(1))}, \nonumber
\end{align}
where $o(1)$ is used for the behavior as $t\rightarrow\infty$. This implies that if $c<\sqrt{\frac{\pi^2}{8 c(\lambda)}}$, then 
\begin{equation} P(\lVert X \rVert_t <c t^{1/3} \mid T>t) \rightarrow 0 \quad \text{as} \quad t\rightarrow \infty. \nonumber
\end{equation}

\vspace{5mm}

\noindent \textbf{Proof of upper bound}

\noindent Let us recall (\ref{eq3}) to start the proof:
\begin{equation} \mathbb{E}^Y\left[\:|R_t(Y+X)|-|R_t(Y)|\:\right]=\mathbb{E}^Y\left[\sup_{s\leq t} (Y_s+X_s)+\sup_{s\leq t} (Y_s-X_s)-2\sup_{s\leq t} (Y_s)\right]. \nonumber
\end{equation}
Let $0<d<1/2$ be a number. Let $c>0$ and $\sigma(X):=\operatorname{argmax}_{{s\in[0,t]}}X_s$. Condition $X$ on the event $\left\{\lVert X \rVert_t\geq c t^d\right\}$. Since $X$ is Brownian, by symmetry, we may suppose without loss of generality that $X_{\sigma(X)}\geq c t^d$. Define 
\[ \tau:=\sup\left\{s\leq \sigma(X):X_{\sigma(X)}-X_s= \frac{ct^d}{2} \right\} \]
so that $0<\tau<\sigma(X)$, and $\tau$ is a random variable depending only on $X$. Note that for $s\in[\tau,\sigma(X)]$, we have $X_s\geq c t^d/2$. Let $\alpha\in(2d,1)$. For each $t\geq 0$, either $|\sigma(X)-\tau|\geq t^\alpha$ or $|\sigma(X)-\tau|\leq t^\alpha$. 

\vspace{0.5cm}

\noindent \underline{Case 1}: Suppose that $|\sigma(X)-\tau|\geq t^\alpha$. Then, for any $t\geq 0$, 
\begin{align} &\mathbb{E}^Y\left[\sup_{s\leq t} (Y_s+X_s)+\sup_{s\leq t} (Y_s-X_s)-2\sup_{s\leq t} (Y_s)\right] \nonumber \\ 
&\:\geq \mathbb{E}^Y\left[\left(\sup_{s\leq t} (Y_s+X_s)+\sup_{s\leq t} (Y_s-X_s)-2\sup_{s\leq t} (Y_s)\right)\mathbbm{1}_{\left\{Y_{\sigma(Y)}\leq c t^d/3\right\}}\mathbbm{1}_{\left\{\sigma(Y)\in[\tau,\sigma(X)]\right\}}\right] \nonumber \\
&\:\geq \mathbb{E}^Y\left[\left(X_{\sigma(Y)}-Y_{\sigma(Y)}\right) \mathbbm{1}_{\left\{Y_{\sigma(Y)}\leq c t^d/3\right\}}\mathbbm{1}_{\left\{\sigma(Y)\in[\tau,\sigma(X)]\right\}}\right] \nonumber \\
&\:\geq \left(c t^d/2-c t^d/3\right)\mathbb{P}^Y\left(Y_{\sigma(Y)}\leq c t^d/3,\:\sigma(Y)\in[\tau,\sigma(X)]\right), \label{eq5}
\end{align}    
where we have used the non-negativity of the expression inside the expectation on the first line of (\ref{eq5}) in passing to the first inequality, and the second supremum on the second line is opened at $s=0$ for a lower bound. It is well known that the joint density of the running maximum $M_t$ and the argument maximum $\sigma_t$ of Brownian motion is given by 
\begin{equation} P_0\left(M_t\in dm, \sigma_t\in du\right)=\frac{m}{\pi}\frac{1}{u^{3/2}\sqrt{t-u}}e^{-m^2/(2u)}dm\: du,\quad m\in[0,\infty),\:u\in(0,t). \label{eq6}
\end{equation}
Since $|\sigma(X)-\tau|\geq t^\alpha$ and $\alpha\in(2d,1)$, integrating (\ref{eq6}) over $u\in[\tau,\sigma(X)]$ and $m\in[0,ct^d/3]$, it is easy to show that for each $\kappa>0$ there exists $c>1$ and $t_1$ such that 
\begin{equation} \mathbb{P}^Y\left(Y_{\sigma(Y)}\leq ct^d/3\:,\sigma(Y)\in[\tau,\sigma(X)]\right)\geq \kappa\:t^{\alpha+2d-2} \label{eq61}
\end{equation}
for all $t\geq t_1$. (Note that as $\kappa$ increases, $c$ increases as well so that we may and do choose $c>1$ in order to bound the first factor on the right-hand side of (\ref{eq5}) from below.) Then, it follows from (\ref{eq5}) that for each $\kappa>0$, there exists $c>0$ and $t_1$ such that conditional on the event $\left\{\lVert X \rVert_t\geq c t^d\right\}$, for all $t\geq t_1$,  
\begin{equation} \mathbb{E}^Y\left[\sup_{s\leq t} (Y_s+X_s)+\sup_{s\leq t} (Y_s-X_s)-2\sup_{s\leq t} (Y_s)\right]\geq  \kappa t^{\alpha+3d-2}. \label{eq7}
\end{equation}

\vspace{0.5cm}

\noindent \underline{Case 2}: Suppose that $|\sigma(X)-\tau|\leq t^\alpha$. Let $\rho=\min\left\{\tau,\sigma(X)-t^{2d}\right\}$ so that $|\sigma(X)-\rho|\geq t^{2d}$. Then, for any $t\geq 0$,
\begin{align} &\mathbb{E}^Y\left[\sup_{s\leq t} (Y_s+X_s)+\sup_{s\leq t} (Y_s-X_s)-2\sup_{s\leq t} (Y_s)\right] \nonumber \\ 
&\:\geq \mathbb{E}^Y\left[\left(\sup_{s\leq t} (Y_s+X_s)+\sup_{s\leq t} (Y_s-X_s)-2\sup_{s\leq t} (Y_s)\right)\mathbbm{1}_{\left\{2 Y_{\sigma(Y)}-Y_{\sigma(X)}-Y_\tau \leq c t^d/3\right\}}\mathbbm{1}_{\left\{\sigma(Y)\in[\rho,\sigma(X)]\right\}}\right] \nonumber \\
&\:\geq \mathbb{E}^Y\left[\left((X_{\sigma(X)}-X_\tau)-(2 Y_{\sigma(Y)}-Y_{\sigma(X)}-Y_\tau)\right) \mathbbm{1}_{\left\{2 Y_{\sigma(Y)}-Y_{\sigma(X)}-Y_\tau \leq c t^d/3 \right\}}\mathbbm{1}_{\left\{\sigma(Y)\in[\rho,\sigma(X)]\right\}}\right] \nonumber \\
&\:\geq \left(c t^d/2-c t^d/3\right)\mathbb{P}^Y\left(2 Y_{\sigma(Y)}-Y_{\sigma(X)}-Y_\tau \leq c t^d/3,\:\sigma(Y)\in[\rho,\sigma(X)]\right), \label{eq8}
\end{align}    
where the first supremum is opened at $s=\sigma(X)$ and the second supremum is opened at $s=\tau$ for a lower bound. Write 
\begin{align} &\mathbb{P}^Y\left(2 Y_{\sigma(Y)}-Y_{\sigma(X)}-Y_\tau \leq c t^d/3,\:\sigma(Y)\in[\rho,\sigma(X)]\right) \nonumber \\
&=\mathbb{P}^Y\left(2 Y_{\sigma(Y)}-Y_{\sigma(X)}-Y_\tau \leq c t^d/3\mid \sigma(Y)\in[\rho,\sigma(X)]\right) \mathbb{P}^Y\left(\sigma(Y)\in[\rho,\sigma(X)]\right). \label{eq9}
\end{align}
If $|\sigma(X)-\tau|\leq t^{2d}$, then the first factor in (\ref{eq9}) is bounded from below by a constant and the second factor by $c t^{2d-1}$ as we shall see below, which means the right-hand side of (\ref{eq9}) is bounded from below by $c t^{2d-1}$. Therefore, suppose that $|\sigma(X)-\tau|\geq t^{2d}$ (so that $\rho=\tau$), which gives a smaller lower bound on (\ref{eq9}) as we show below. Since now $t^{2d}\leq|\sigma(X)-\tau|\leq t^\alpha$, for a lower bound on (\ref{eq9}) for large $t$, we may take $|\sigma(X)-\tau|= t^\alpha$. (It is clear that the first factor on the right-hand side of (\ref{eq9}) decreases as $|\sigma(X)-\tau|$ increases, whereas the second factor increases, but the effect on the first factor is dominating for large $t$ as the product of lower bounds in (\ref{eq11}) and (\ref{eq19}) reveals below.) 

It is well known that the argument maximum $\sigma_t$ of Brownian motion has the arcsine distribution given by 
\begin{equation} P\left(\sigma_t\in du\right)=\frac{1}{\pi\sqrt{u(t-u)}}du,\quad u\in(0,t). \nonumber
\end{equation} 
Since we have taken $|\sigma(X)-\tau|=t^\alpha$, it follows that for any $t>0$,
\begin{equation} \mathbb{P}^Y\left(\sigma(Y)\in[\tau,\sigma(X)]\right)\geq \frac{2}{\pi t}\: t^\alpha=\frac{2}{\pi}\: t^{\alpha-1}, \label{eq11}
\end{equation}
where $2/(\pi t)$ is the minimum of the arcsine distribution on $(0,t)$.

Next, we wish to find a lower bound for 
\begin{equation} \mathbb{P}^Y\left(2 Y_{\sigma(Y)}-Y_{\sigma(X)}-Y_\tau \leq c t^d/3\mid \sigma(Y)=r \right) \label{kule}
\end{equation}
uniformly in $r\in[\tau,\sigma(X)]$. It is easy to see from the time-reversal symmetry of Brownian motion that since $|\sigma(X)-\tau|= t^\alpha$, for $r\in[\tau,\sigma(X)]$ and any $t>0$,
\begin{equation} \mathbb{P}^Y\left(2 Y_{\sigma(Y)}-Y_{\sigma(X)}-Y_\tau \leq c t^d/3\mid \sigma(Y)=r \right) \geq \inf_{r\leq t-t^\alpha}\mathbb{P}^Y\left(Y_{\sigma(Y)}-Y_{\sigma(Y)+t^{\alpha}} \leq c t^d/6 \mid \sigma(Y)=r\right). \label{eq12}
\end{equation}
Observe that conditional on $\left\{\sigma(Y)=r\right\}$, the process $(Y_{\sigma(Y)+s}-Y_{\sigma(Y)})_{s\in[0,t-r]}$ is a Brownian motion starting at $0$, and conditioned to avoid $[0,\infty)$. Let $\tau_{A}=\inf\left\{s>0:Y_s\in A\right\}$ be the first hitting time of the set $A\subset\mathbb{R}$ by $Y$. Set $\tau_0:=\tau_{(-\infty,0]}$. It then follows by symmetry that for any $t>0$, (\ref{kule}) is bounded from below by
\begin{equation} \inf_{r\leq t-t^\alpha}\mathbb{P}^Y\left(Y_{t^\alpha}\leq c t^d/6\mid \:\tau_0>t-r\right). \label{eq13}
\end{equation}
For $t>1$, to avoid working with an event that has zero probability (the event $\left\{\tau_0>t-r\right\}$ is as such), suppose that $Y$ is started at $1$ instead of $0$. Note that $Y$ started at any $0<x<c t^d/6$ instead of $0$ can only decrease the probability in (\ref{eq13}). Let $\mathbb{P}^Y_y$ be the law of $Y$ started at $y$. For $0<s<t$, the transition probability density for Brownian motion started at $x>0$ at time $0$ and arriving at $y>0$ at time $s$, and that is conditioned to stay positive up to time $t$, is given by (see for instance \cite[ex.1.14(ii)]{K2016}) 
\begin{equation} p_t(s,y\mid x):=\frac{1}{\sqrt{2\pi s}}\left[e^{-(x-y)^2/(2s)}-e^{-(x+y)^2/(2s)}\right]\frac{P_y(\tau_0>t-s)}{P_x(\tau_0>t)}, \quad y\in(0,\infty), \nonumber
\end{equation}
where $P_x$ is, as introduced before, the law of Brownian motion started at $x$. It follows that for any $t>0$ and $r\leq t-t^{\alpha}$,
\begin{align} &\mathbb{P}^Y\left(Y_{t^\alpha}\leq c t^d/6\mid \:\tau_0>t-r\right) \nonumber \\
& \geq \frac{1}{\sqrt{2\pi t^\alpha}}\frac{\int_0^{c t^d/6}\left[e^{-(1-y)^2/(2t^\alpha)}-e^{-(1+y)^2/(2t^\alpha)}\right] \mathbb{P}^Y_y\left(\tau_0>t-r-t^\alpha\right) dy}{\mathbb{P}^Y_1\left(\tau_0>t-r\right)} \nonumber \\
& \geq \frac{1}{\sqrt{2\pi t^\alpha}}\frac{\mathbb{P}^Y_{c t^d/7}\left(\tau_0>t-r-t^\alpha\right)}{\mathbb{P}^Y_1\left(\tau_0>t-r\right)}\int_{c t^d/7}^{c t^d/6}\left[e^{-(1-y)^2/(2t^\alpha)}-e^{-(1+y)^2/(2t^\alpha)}\right] dy, \label{eq15}
\end{align}
where in passing to the last inequality, firstly the lower limit of integration was shifted from $0$ to $c t^d/7$ since the integrand is positive, and then $\mathbb{P}^Y_y\left(\tau_0>t-r-t^\alpha\right)$ was taken outside the integral by setting $y=c t^d/7$ for a lower bound. Recall that the probability density of first hitting time of zero for a Brownian motion started at $x>0$ is given by 
\begin{equation} P_x\left(\tau_0\in du\right)=\frac{xe^{-x^2/(2u)}}{\sqrt{2\pi u^3}}du, \quad u\in(0,\infty). \nonumber
\end{equation}
It follows that for all large $t$ and $r\leq t-t^{\alpha}$,
\begin{align} \frac{\mathbb{P}^Y_{c t^d/7}\left(\tau_0>t-r-t^\alpha\right)}{\mathbb{P}^Y_1\left(\tau_0>t-r\right)}&=\frac{c t^d}{7}\frac{\int_{t-r-t^\alpha}^\infty \frac{1}{u^{3/2}}\exp[-(c t^d/7)^2/(2u)]du}{\int_{t-r}^\infty \frac{1}{u^{3/2}}\exp[-1/(2u)]du} \nonumber \\
&\geq \frac{c t^d}{7}\frac{\int_{t-r}^\infty \frac{1}{u^{3/2}}\exp[-(c t^d/7)^2/(2u)]du}{\int_{t-r}^\infty \frac{1}{u^{3/2}}\exp[-1/(2u)]du}. \nonumber \\
&\geq \frac{c t^d}{14}. \label{eq17}
\end{align}
In passing to the last inequality, we have used that for all large $t$, $\exp[-(c t^d/7)^2/(2u)]\geq 1/2$ uniformly in $u\in[t-r,\infty)$ since $u\geq t-r\geq t^\alpha$ and $\alpha>2d$, and that $\exp[-1/(2u)]\leq 1$. Furthermore, it is easy to see that for each $\kappa>0$ there exists $c>1$ and $t_2$ such that  
\begin{equation} \int_{c t^d/7}^{c t^d/6}\left[e^{-(1-y)^2/(2t^\alpha)}-e^{-(1+y)^2/(2t^\alpha)}\right] dy \geq \kappa t^{2d-\alpha}, \nonumber
\end{equation}
for all $t\geq t_2$, and hence by (\ref{eq15}) and (\ref{eq17}) that 
\begin{equation} \mathbb{P}^Y\left(Y_{t^\alpha}\leq c t^d/6\mid \:\tau_0>t-r\right)\geq \kappa t^{3d-3\alpha/2} \label{eq19}
\end{equation}
for all $t\geq t_2$ uniformly in $r\leq t-t^\alpha$. Combining this with (\ref{eq8})-(\ref{eq11}) and (\ref{eq12}), it follows that for each $\kappa>0$, there exists $c>0$ and $t_2$ such that conditional on the event $\left\{\lVert X \rVert_t\geq c t^d\right\}$, for all $t\geq t_2$,  
\begin{equation} \mathbb{E}^Y\left[\sup_{s\leq t} (Y_s+X_s)+\sup_{s\leq t} (Y_s-X_s)-2\sup_{s\leq t} (Y_s)\right]\geq  \kappa t^{4d-1-\alpha/2}. \label{eq20}
\end{equation}

Recall that the choice of $\alpha\in(2d,1)$ was arbitrary. Optimizing the lower bounds in (\ref{eq7}) and (\ref{eq20}) over $\alpha$ gives $\alpha=2(d+1)/3$, which in turn gives for both cases $|\sigma(X)-\tau|\geq t^\alpha$ and $|\sigma(X)-\tau|\leq t^\alpha$ the following result: for each $\kappa>0$, there exists $c>0$ and $t_0$ such that conditional on the event $\left\{\lVert X \rVert_t\geq c t^d\right\}$, for all $t\geq t_0$,
\begin{equation} \mathbb{E}^Y\left[\sup_{s\leq t} (Y_s+X_s)+\sup_{s\leq t} (Y_s-X_s)-2\sup_{s\leq t} (Y_s)\right]\geq  \kappa t^{\frac{11d}{3}-\frac{4}{3}}. \label{eq21}
\end{equation}
Since $\kappa$ in (\ref{eq21}) can be arbitrarily large, in view of (\ref{eq1}) and (\ref{eq3}), to prove the upper bound in (\ref{eq01}), it suffices that $\frac{11d}{3}-\frac{4}{3}=\frac{1}{3}$, which yields $d=5/11$. This completes the proof.

\section{Proof of Theorem~\ref{teo3} and Theorem~\ref{teo2} }\label{s5}

\subsection{Proof of Theorem~\ref{teo3}}

Let $0<\varepsilon\leq 1$ and $0<\kappa<1$ be fixed. Suppose that $X$ stays inside $(\kappa\lVert X \rVert_t,\lVert X \rVert_t)$ for a total time of length $\varepsilon t$ so that $S:=\left\{s\in[0,t]:X_s\in(\kappa\lVert X \rVert_t,\lVert X \rVert_t)\right\}$ has Lebesgue measure $\geq \varepsilon t$. Let $c_3>1$ be a constant which will depend on $\kappa$ and $\varepsilon$. Then, conditional on $|S|\geq \varepsilon t$ and $\lVert X \rVert_t \geq c_3\,t^{4/9}$,   
\begin{align} &\mathbb{E}^Y\left[\sup_{s\leq t} (Y_s+X_s)+\sup_{s\leq t} (Y_s-X_s)-2\sup_{s\leq t} (Y_s)\right] \nonumber \\ 
&\:\geq \mathbb{E}^Y\left[\left(\sup_{s\leq t} (Y_s+X_s)+\sup_{s\leq t} (Y_s-X_s)-2\sup_{s\leq t} (Y_s)\right)\mathbbm{1}_{\left\{Y_{\sigma(Y)}\leq \kappa c_3 t^{4/9}/2\right\}}\mathbbm{1}_{\left\{\sigma(Y)\in S\right\}}\right] \nonumber \\
&\:\geq \mathbb{E}^Y\left[\left(X_{\sigma(Y)}-Y_{\sigma(Y)}\right) \mathbbm{1}_{\left\{Y_{\sigma(Y)}\leq \kappa c_3 t^{4/9}/2\right\}}\mathbbm{1}_{\left\{\sigma(Y)\in S\right\}}\right] \nonumber \\
&\:\geq \left(\kappa c_3 t^{4/9}-\kappa c_3 t^{4/9}/2\right)\mathbb{P}^Y\left(Y_{\sigma(Y)}\leq \kappa c_3 t^{4/9}/2,\:\sigma(Y)\in S\right) \label{eq23}
\end{align}    
for any $t\geq 0$. One can show similarly to the argument that follows (\ref{eq5}) that for every $c_4$, there exists $c_3>1$ such that
\begin{equation} \mathbb{P}^Y\left(Y_{\sigma(Y)}\leq \kappa c_3 t^{4/9}/2\:,\sigma(Y)\in S\right)\geq c_4\:t^{-1/9} \nonumber
\end{equation} 
for all large $t$. This completes the first part of the proof in view of (\ref{eq1}), (\ref{eq3}) and (\ref{eq23}). (Note that the argument here is the extension of case $1$ in the proof of the upper bound of Theorem~\ref{teo1} to $\alpha=1$.)

Now let $A_t$ be the event that $X$ traverses the interval $(\kappa\lVert X \rVert_t,\lVert X \rVert_t)$ at least diffusively fast so that there exists $k>0$ and times $\tau_1,\tau_2\in[0,t]$ with $X_{\tau_1}=\kappa\lVert X \rVert_t$ and $X_{\tau_2}=\lVert X \rVert_t$ such that $|\tau_1-\tau_2|\leq k((1-\kappa)\lVert X \rVert_t)^2$. Let $c_3>1$ be a constant which will depend on $\kappa$ and $k$. Then, conditional on $A_t$ and $\lVert X \rVert_t \geq c_3\,t^{4/9}$,   
\begin{align} &\mathbb{E}^Y\left[\sup_{s\leq t} (Y_s+X_s)+\sup_{s\leq t} (Y_s-X_s)-2\sup_{s\leq t} (Y_s)\right] \nonumber \\ 
&\:\geq \mathbb{E}^Y\left[\left(\sup_{s\leq t} (Y_s+X_s)+\sup_{s\leq t} (Y_s-X_s)-2\sup_{s\leq t} (Y_s)\right)\mathbbm{1}_{\left\{2Y_{\sigma(Y)}-Y_{\tau_1}-Y_{\tau_2}\leq (1-\kappa)\lVert X \rVert_t/2\right\}}\mathbbm{1}_{\left\{\sigma(Y)\in U\right\}}\right] \nonumber \\
&\:\geq \mathbb{E}^Y\left[\left((X_{\tau_2}-X_{\tau_1})-(2Y_{\sigma_Y}-Y_{\tau_1}-Y_{\tau_2})\right) \mathbbm{1}_{\left\{2Y_{\sigma(Y)}-Y_{\tau_1}-Y_{\tau_2}\leq (1-\kappa)\lVert X \rVert_t/2\right\}}\mathbbm{1}_{\left\{\sigma(Y)\in U\right\}}\right] \nonumber \\
&\:\geq \left((1-\kappa)c_3 t^{4/9}/2\right)\mathbb{P}^Y\left(2Y_{\sigma(Y)}-Y_{\tau_1}-Y_{\tau_2}\leq (1-\kappa)\lVert X \rVert_t/2 ,\:\sigma(Y)\in U\right) \label{eq25}
\end{align} 
for any $t\geq 0$, where the interval $U$ is chosen such that $[\tau_1,\tau_2]\subseteq U \subseteq [0,t]$ and $|U|\geq\min\left\{(\lVert X \rVert_t)^2,t\right\}$. Since $|U|\geq c_3 t^{8/9}$, one can show similarly to the argument that follows (\ref{eq9}) that for every $c_4$, there exists $c_3>1$ such that conditional on $A_t$ and $\lVert X \rVert_t \geq c_3\,t^{4/9}$,
\begin{equation} \mathbb{P}^Y\left(2Y_{\sigma(Y)}-Y_{\tau_1}-Y_{\tau_2}\leq (1-\kappa) \lVert X \rVert_t/2\:,\sigma(Y)\in U\right)\geq c_4\:t^{-1/9} \nonumber
\end{equation} 
for all large $t$. This completes the proof in view of (\ref{eq1}), (\ref{eq3}) and (\ref{eq25}).

\subsection{Proof of Theorem~\ref{teo2}}

Let $\varepsilon>0$ be fixed and $f:\mathbb{R}_+\to\mathbb{R}_+$ be a function such that $f(t)\rightarrow\infty$ and $f(t)=o(t^{1/3})$ as $t\rightarrow\infty$. Let $I_n(t)$, $n=1,2,\ldots,\left\lceil \varepsilon t^{1/3}/f(t)\right\rceil$ be pairwise disjoint intervals in $[0,t]$, each of length $|I_n|= t^{2/3}f(t)$. In the rest of the proof, $n$ runs from $1$ to $\left\lceil \varepsilon t^{1/3}/f(t)\right\rceil$, and phrases such as `for all $n$' will mean for $1\leq n \leq \left\lceil \varepsilon t^{1/3}/f(t)\right\rceil$. Let $R_{I_n}:=\left\{X_s:s\in I_n\right\}$. Let $B_t$ be the event that there exists $k>0$ such that $|R_{I_n}(X)|\geq kt^{1/3}\sqrt{f(t)}$ for all $n$, so that conditioned on $B_t$, $X$ is diffusive on $\geq \varepsilon t^{1/3}/f(t)$ many pairwise disjoint intervals in $[0,t]$ of length $t^{2/3}f(t)$ each.

Let $S=\cup_{n} I_n$. Then, conditioned on $B_t$, for any $t\geq 0$, 
\begin{align} &\mathbb{E}^Y\left[\sup_{s\leq t} (Y_s+X_s)+\sup_{s\leq t} (Y_s-X_s)-2\sup_{s\leq t} (Y_s)\right] \nonumber \\ 
&\:\geq \mathbb{E}^Y\left[\left(\sup_{s\leq t} (Y_s+X_s)+\sup_{s\leq t} (Y_s-X_s)-2\sup_{s\leq t} (Y_s)\right)\mathbbm{1}_{\left\{2Y_{\sigma(Y)}-Y_u-Y_v\leq kt^{1/3}\sqrt{f(t)}/2\right\}}\mathbbm{1}_{\left\{\sigma(Y)\in S\right\}}\right] \nonumber \\
&\:\geq \mathbb{E}^Y\left[\left(X_u-X_v\right)-\left(2Y_{\sigma(Y)}-Y_u-Y_v\right)\mathbbm{1}_{\left\{2Y_{\sigma(Y)}-Y_u-Y_v\leq kt^{1/3}\sqrt{f(t)}/2\right\}}\mathbbm{1}_{\left\{\sigma(Y)\in S\right\}}\right] \nonumber \\
&\:\geq \left(kt^{1/3}\sqrt{f(t)}-kt^{1/3}\sqrt{f(t)}/2\right)\mathbb{P}^Y\left(2Y_{\sigma(Y)}-Y_u-Y_v\leq kt^{1/3}\sqrt{f(t)}/2,\:\sigma(Y)\in S\right), \label{eq22}
\end{align}    
where $u$ and $v$ are chosen from the interval $I_j$ into which $\sigma(Y)$ falls, in such a way that $X_u-X_v\geq kt^{1/3}\sqrt{f(t)}$. Note that since we condition on $B_t$, by definition of $B_t$, it is possible to find such a pair $u,v$ in each $I_n$. Since $\max_{w\in\left\{u,v\right\}}\left\{|\sigma(Y)-w|\right\}\leq t^{2/3}f(t)$ and $|S|\geq \varepsilon t$, the second factor on the right-hand side of (\ref{eq22}) can be bounded from below by a constant $c(\varepsilon,k)$ for all large $t$ by using a method similar to the one used in the proof of Theorem~\ref{teo1} starting with (\ref{eq9}). This completes the proof in view of (\ref{eq1}) and (\ref{eq3}).

\section{Open problems}\label{s6}

We conclude by giving several open problems related to our model.

1. Sharp upper bound on maximal displacement in $d=1$: Here, we do not claim that our upper bound in (\ref{eq01}) is sharp. Further work is needed to either show that the exponent $5/11$ is sharp or find a better upper bound. We note that in \cite{ADS2016}, it was conjectured that the fluctuations of $X$ are on the scale of $t^{1/3}$, which would mean that the exponent $5/11$ in (\ref{eq01}) could be lowered to $1/3$.   

2. Maximal displacement in higher dimensions: The current work, in particular Theorem~\ref{teo1}, is for $d=1$. The maximal displacement of $X$ from origin in $d\geq 2$ conditioned on survival stands as an open problem. In $d=1$, we are able to express $|W_X(t)|$, that is, the volume of the Wiener sausage `perturbed by' $X$, in terms of the running maximum and running minimum of Brownian motions, on which the entire analysis is based. In $d\geq 2$, this nice connection to running extrema is lost, therefore a different approach is needed.

3. Other kinds of optimal survival strategies: Finer questions could be asked about the optimal survival strategy followed by the Brownian particle $X$. For instance, what proportions of the time interval $[0,t]$ does $X$ spend `near' origin and `far away' from origin? \par Strategies involving the trap field are also of interest. Our system is composed of $X$ and the trap field, and we have considered the strategies involving $X$ only. A natural question is: Does the trap field leave out space-time clearings (trap-free regions) in all or part of $[0,t]$ with overwhelming probability, given that $X$ avoids traps up to time $t$? Recall that in Section~\ref{s3}, in order to find a lower bound for $P(T>t)$, we have used a strategy where the trap field avoids $B(0,t^{1/3}+a)$ throughout $[0,t]$. Is this strategy optimal or does clearing a ball with radius $o(t^{1/3})$ in some or all of $[0,t]$ a better strategy for survival (possibly coupled with a strategy followed by $X$)?  
    
4. Survival asymptotics in $d\geq 3$: To the best of our knowledge, an exact value for 
	\begin{equation} -\underset{t\rightarrow\infty}{\lim}\frac{1}{t}\log(\mathbb{E}\times\mathbb{P}^X)(T>t) \quad \text{in} \quad d\geq 3 \nonumber
	\end{equation}
has not been found. One can show using a subadditivity argument that this Lyapunov exponent exists and is positive, and can bound it from below using Jensen's inequality and from above using Pascal's principle (see \cite{PS2012}). We conjecture that the upper bound coming from Pascal's principle coincides with the actual value of the exponent, however we are unable to prove it. We note that the analogous question in the discrete version of the current problem was also left open in \cite{R2012}.

\bibliographystyle{plain}

\end{document}